# A new probability model with support on unit interval: Structural properties, regression of bounded response and applications

Subrata Chakraborty[a]*, Seng Huat Ong[b] and C. M. Ng[c]

[a]*Department of Statistics, Dibrugarh University, Dibrugarh, India;*
[b]*Department of Actuarial Science and Applied Statistics, UCSI University, Kuala Lumpur, Malaysia;*
[c]*Institute of Mathematical Sciences, University of Malaya, Malaysia*

**email: subrata_arya@yahoo.co.in*

**Abstract**

A new distribution on (0, 1), generalized Log-Lindley distribution, is proposed by extending the Log-Lindley distribution. This new distribution is shown to be a weighted Log-Lindley distribution. Important probabilistic and statistical properties have been derived. An interesting characterization of the weighted distribution in terms of Kullback-Liebler distance and weighted entropy has also been obtained. A useful result in insurance for the distorted premium principal is presented and verified with numerical calculations. New regression models for bounded responses based on this distribution and their application is illustrated by considering modeling a real life data on risk management and another data set on outpatient health expenditure in comparison with beta regression and Log-Lindley regression models. Much better fits for both data sets justify the relevance of the new distribution in statistical modeling and analysis. Furthermore this generalization, apart from adding flexibility for modelling, retains the compactness and tractability of statistical quantities required for statistical analysis, which is a feature of the Log-Lindley distribution. Thus, the generalized Log-Lindley distribution should be a useful addition to statistical models for practitioners.

***Keywords***: computer sampling; distortion function; log-concavity; premium loading; regression; stochastic ordering; weighted entropy

***Subject classification codes***: 60EXX, 62E15, 62F10, 62J12

**(Version: 4th February 2020)**

## 1. Introduction

In statistical research many attempts have been made to introduce alternatives to the classical beta distribution (see [1] for details). But most of these distributions involved special functions in their formulations except for the Kumaraswamy distribution [2]. Recently, a new probability density function (pdf) with bounded support in (0,1) was introduced by Gomez-Deniz et al. [3], as an alternative to the classical beta distribution by suitable transformation of a particular case of the generalized Lindley distribution proposed by Zakerzadeh and Dolati [4]. The new distribution called Log-Lindley (LL) distribution has compact expressions for the moments as well as the cdf. Gomez-Deniz et al. [3] studied its important properties relevant to the insurance and inventory management applications. In the application to insurance premium loading, Gomez-Deniz et al. [3] showed that the risk-adjusted (distorted) premium based on the LL distribution falls between the net premium and the dual-power premium. The LL model is also shown to be appropriate regression model to model bounded responses as an alternative to the beta regression model. The LL distribution of Gomez-Deniz et al. [3] is probably the latest in a series of proposals as alternatives to the classical beta distribution. Unlike its predecessors the main attraction of the LL distribution is that it has nice compact forms for the pdf, cdf and moments, which do not involve any special functions. Recently, Jodra and Jimenez-Gamero [5] derived the quantile function of the LL distribution in terms of the Lambert W function and proposed this as a method to sample from the LL distribution.

The pdf, cdf and moments of the LL distribution of Gomez-Deniz et al. [3] are given by

Pdf: $$f(x;\theta,\lambda) = \frac{\theta^2}{1+\lambda\theta}(\lambda - \log x)\, x^{\theta-1}, 0 < x < 1,\ \lambda \geq 0,\ \theta > 0 \qquad (1)$$

Cdf: $$F(x;\theta,\lambda) = \frac{x^{\theta}[1+\theta(\lambda - \log x)]}{1+\lambda\theta}$$

Moments: $E(X^r;\theta,\lambda)=\dfrac{\theta^2}{1+\lambda\theta}\dfrac{1+\lambda(r+\theta)}{(r+\theta)^2}, r=\cdots,-2,-1,1,2,\cdots$

In this paper we extend the LL distribution to provide another distribution on (0, 1). A probability distribution defined on (0, 1) is useful in many ways, for instance, in insurance as a distortion function in the definition of a premium principle. The motivation to propose this extension is that apart from adding flexibility to the potentially very useful LL distribution as an alternative to the beta distribution, it also retains the compactness and tractability of expressions for the pdf, cdf and moments which facilitate its applications in statistical analysis. Furthermore, a number of important probabilistic properties have been derived and the proposed distribution is shown to be a weighted LL distribution. We also justified the proposed extension of the LL distribution by showing improvement in fit to the same data used by Gomez-Deniz et al. [3] to support the application of their LL model. This data set was originally used in Schmit and Roth [6] and is about the cost effectiveness of risk management (measured in percentages) in relation to exposure to certain property and casualty losses, adjusted by several other variables such as size of assets and industry risk. Data were collected from a questionnaire survey of 374 risk managers (73 observations) of large US based organizations. Further details are given in Section 6.1. A more recent data set from the Medical Expenditure Panel Survey (MEPS), conducted by the U.S. Agency of Health Research and Quality has also been fitted by the proposed extension of the LL distribution. MEPS is a probability survey that gives nationally representative estimates of health care use, expenditures, sources of payment, and insurance coverage for the U.S. civilian population. Detailed information on individuals of each medical care episode by type of services is collected and this allows the development of models of health care utilization to predict future expenditures. With the added flexibility of the extended LL model, a much better fit than the beta and Log-Lindley distributions is also obtained.

The paper has been organized as follows. Section 2 presents the generalized LL distribution and various properties like moments, mode, Shannon entropy, stochastic ordering and convexity and log-concavity. The distribution is derived as a weighted LL distribution. Computer sampling from the generalized LL distribution is also examined in Section 2 and as a special case we propose a faster method of computer generation of LL samples. An interesting characterization of the weighted distribution in terms of Kullback-Liebler distance and weighted entropy has been obtained. An application to insurance premium loading is considered in Section 3 where it is shown that the premium based on the GLL distribution lies between the proportional hazard premium [7,8] and maximal premium. Section 4 gives a useful re-parameterization of the generalized Log-Lindley distribution. Section 5 gives the score equations for maximum likelihood estimation and simple expressions for the elements of the information matrix. Illustration of non-regression and regression data modeling with the risk management data set in [6], and used by Gomez-Deniz et al. [3] for the LL distribution, is examined in Section 6. A concluding discussion ends the paper.

## 2. Generalized Log-Lindley distribution $LL_p(\theta, \lambda)$

The LL distribution of Gomez-Deniz et al. [3] was obtained from a two-parameter particular case of the three-parameter generalized Lindley (GL) distribution of Zakerzadeh and Dolati [4] which has the pdf

$$f(x) = \frac{\theta^2 (\theta x)^{\alpha-1} (\alpha + \lambda x) e^{-\theta x}}{(\lambda + \theta) \Gamma(\alpha + 1)}, 0 < x < 1,\ \lambda \geq 0,\ \theta > 0.$$

Gomez-Deniz et al. [3] first substituted $\alpha = 1$ and then employed the transformation $X = \log(1/Y)$ to obtain their LL distribution in (1).

### *2.1 Definition and derivation as a weighted log-Lindley distribution*

In this section, we consider the same transformation of the GL distribution but without setting

$\alpha = 1$ to derive an extension called the **generalized Log-Lindley** (GLL) distribution. This proposed three-parameter distribution can be seen as an extension of the LL distribution of Gomez-Deniz et al. (2014) with an additional parameter '*p*'. We denote the new distribution by $LL_p(\theta, \lambda)$. By setting $p = \alpha - 1$, the pdf and cdf of the proposed distribution are given respectively by

$$f(x;\theta,\lambda,p) = \frac{\theta^{2+p}}{\Gamma(1+p)[1+p+\lambda\theta]}(-\log x)^p(\lambda - \log x)\,x^{\theta-1}, \;\; 0<x<1,\; \lambda \ge 0,\; \theta > 0,\; p \ge 0 \tag{2}$$

and $$F(x;\theta,\lambda,p) = \frac{[x^{\theta} + (1+p+\theta\lambda)Ei(-p,-\theta\log x)](-\theta\log x)^{1+p}}{\Gamma(1+p)[1+p+\lambda\theta]} \tag{3}$$

where $Ei(n,z) = \int_1^{\infty} e^{-zt}/t^n dt = z^{n-1}\Gamma(1-n,z)$, $z > 0$, is the generalized exponential integral and $\Gamma(1-n,z)$ is the upper incomplete gamma function.

In particular, the recurrence relations $n\,Ei(n+1,z) = e^{-z} - zEi(n,z)$ and $\frac{d}{dz}Ei(n,z) = -zEi(n-1,z)$ can be utilized in conjunction with $Ei(0,z) = e^{-z}/z$ to obtain expressions for higher integral values on *n*. Moreover, as $z \to \infty$, $Ei(n,z) = e^{-z}/z$ (see http:// mathworld. wolfram.com/ En-Function.html for more results).

Alternatively, we can write the cdf in (3) in terms of upper incomplete gamma function as

$$F(x;\theta,\lambda,p) = \frac{x^{\theta}(-\theta\log x)^{1+p} + (1+p+\theta\lambda)\Gamma(1+p,-\theta\log x)}{\Gamma(1+p)[1+p+\lambda\theta]}.$$

Further using the result $\Gamma(p+1,y) = p!\,e^{-y}\sum_{k=0}^{p} y^k/k!$ for integer *p*, the cdf can be written in a compact form devoid of any special function as

$$F(x;\theta,\lambda,p)=\frac{x^{\theta}(-\theta\log x)^{1+p}+(1+p+\theta\lambda)p!\,x^{\theta}\sum_{k=0}^{p}(-\theta\log x)^{k}/k!}{(1+p+\lambda\theta)\,p!}.$$

Throughout the article $LL_p(\theta,\lambda)$ and **GLL** have been used synonymously.

We now formulate the generalized log-Lindley distribution as a weighted log-Lindley distribution. The general concept of weighted probability distributions was formalized by Rao [9] and has found applications in diverse areas. A weighted distribution is a result of modifying the original distribution by a weight function. This arises from the fact that the process of recording observations is done without equal probability for each observation. Let $X$ be a random variable with pdf $f(x)$ and $w(x)$ be a weight function. The weighted pdf $f_W(x)$ is given by

$$f_W(x)=\frac{w(x)f(x)}{E\left[w(X)\right]}.$$

Let $f(x)$ be the log-Lindley pdf and $w(x)=(-\log x)^{p}$. We find that

$$E\left[w(X)\right]=\frac{\theta^{p}(1+\lambda\theta)}{\Gamma(1+p)(1+p+\lambda\theta)}.$$

Then

$$f_W(x)=\frac{\theta^{p}(1+\lambda\theta)}{\Gamma(1+p)(1+p+\lambda\theta)}(-\log x)^{p}\frac{\theta^{2}}{(1+\lambda\theta)}(\lambda-\log x)x^{\theta-1}.$$

This is the GLL pdf given by (2). □

***2.2 Important particular cases of $LL_p(\theta,\lambda)$ and shape of $LL_p(\theta,\lambda)$ distributions***

For

**(1)** $p=0$, we get back the $LL$ distribution of Gomez-Deniz et al. [3] given by (1).

**(2)** $p=0,\theta=1$, as $\lambda\to\infty$, $LL_p(\theta,\lambda)\to Uniform(0,1)$.

**(3)** $p=1$, we get new $LL_1(\theta,\lambda)$ distribution as:

Pdf: $f(x;\theta,\lambda)=\dfrac{\theta^3}{[2+\lambda\theta]}(-\log x)(\lambda-\log x)\,x^{\theta-1}, \quad 0<x<1,\ \lambda\geq 0,\ \theta>0$

Cdf: $F(x;\theta,\lambda)=\dfrac{x^\theta[2+\theta\lambda-\theta\log x(2+\theta\lambda-\theta\log x)]}{2+\lambda\theta}$

Moment: $E(X^r;\theta,\lambda)=\dfrac{\theta^3[(r+\theta)\lambda+2]}{(r+\theta)^3[2+\lambda\theta]}, \quad r=\cdots,-2,-1,1,2,\cdots$

**(4)** $p=2$, we get new $LL_2(\theta,\lambda)$ distribution as:

Pdf: $f(x;\theta,\lambda)=\dfrac{\theta^4}{6+2\lambda\theta}(-\log x)^2(\lambda-\log x)\,x^{\theta-1}, \quad 0<x<1,\ \lambda\geq 0,\ \theta>0$

Cdf: $F(x;\theta,\lambda)=\dfrac{x^\theta[6+2\theta\lambda-\theta\log x\{6+2\theta\lambda-\theta\log x(3+\theta\lambda-\theta\log x)\}]}{6+2\theta\lambda}$

Moment: $E(X^r;\theta,\lambda)=\dfrac{\theta^4[2\lambda(r+\theta)+6]}{2(r+\theta)^4(3+\lambda\theta)}, \quad r=\cdots,-2,-1,1,2,\cdots$

Here we have illustrated the pdf of $LL_p(\theta,\lambda)$ distributions for different choices of parameters $\theta,\lambda$ and index $p$ to study their shapes. It is observed that the distribution can be positively or negatively skewed and symmetrical. The pdf can also be increasing or decreasing.

**Figure 1 here**

### *2.3 Moments and Mode*

The moments of $LL_p(\theta,\lambda)$ distribution are given by

$$E(X^r;\theta,\lambda,p)==\frac{\theta^{2+p}((r+\theta)\lambda+(1+p))}{(r+\theta)^{2+p}[1+p+\lambda\theta]}, r=\cdots,-2,-1,1,2,\cdots$$

In particular, the mean is

$$E(X;\theta,\lambda,p)=\left(\frac{\theta}{1+\theta}\right)^{2+p}\frac{1+p+\lambda(1+\theta)}{1+p+\lambda\theta}=\mu\,. \qquad (4)$$

The mode plays an important role in the usefulness of a distribution. For the $LL_p(\theta,\lambda)$ distribution the mode occurs at

$$\exp[\{1+p+\lambda(1-\theta)\pm\sqrt{4p\lambda(\theta-1)+(1+p+\lambda(1-\theta))^2}\}/2(1-\theta)], \text{ for } \theta\neq 1.$$

### *2.4 Random sampling from the GLL Distribution*

Jodra and Jimenez-Gamero [5] derived the quantile function of the LL distribution in terms of the Lambert W function. This gives a direct method of generating LL samples on the computer which only requires the computation of the Lambert W function. Computation of special functions usually require relatively more computational effort compared to elementary functions. This clearly impacts on the speed of the computer sampling which is crucial in a massive Monte Carlo simulation experiment. We give a quick method of generating samples from the LL distribution which is a particular case of the method for computer generation of the generalized Log-Lindley distribution presented below.

Zakerzadeh and Dolati [4] expressed the pdf of their generalized Lindley distribution as a finite mixture of two gamma pdf's:

$$f(x)=\frac{\theta}{\lambda+\theta}f_g(x;\alpha,\theta)+\frac{\lambda}{\lambda+\theta}f_g(x;\alpha+1,\theta)$$

where $f_g(x;\alpha,\theta)=\theta^{\alpha}x^{\alpha-1}e^{-\theta x}/\Gamma(\alpha)$.

Based on this finite mixture representation and the transformation $Z=-\log(X)$, we propose the following method to generate samples from the GLL distribution. Let $G(\alpha,\theta)$ denote the gamma random variable with parameters $(\alpha,\theta)$ and $U$ be the uniform random variable over $(0,1)$.

***Algorithm 1 (Generation of GLL samples)***

Specify parameters $(\lambda,\alpha,\theta)$.

**(1)** Generate $U$.

**(2)** If $U \le \frac{\theta}{\lambda+\theta}$, generate $X = G(\alpha,\theta)$; otherwise generate $X = G(\alpha+1,\theta)$.

**(3)** Accept $Z = -\log(X)$ as a GLL random variate (transform $X$ ).

For a sample of size $N$, repeat Steps 1, 2 and 3 $N$ times.

There are efficient methods for computer sampling from the gamma distribution. For example, Cheng [10] gave an efficient algorithm based on the acceptance-rejection method.

If $\alpha = 1$, ***Algorithm 1*** provides a method of generating LL random samples. We note that when $\alpha = 1$, $G(1,\theta)$ is the exponential random variable and the quantile function is $Q(u) = -\theta^{-1}\log(1-u)$, $0 < u < 1$. Therefore, generation of an exponential random variate is very quick requiring only evaluation of the log function. Since $1-U$ is also uniform on $(0,1)$, an exponential random variate is given by $X = -\theta^{-1}\log(U)$, and this saves one arithmetic operation. The random variable $X = G(2,\theta) = G(1,\theta) + G(1,\theta)$ may be taken as a sum of two independent exponential random variables. We have the following algorithm for the generation of LL samples.

***Algorithm 2 (Generation of Log-Lindley samples)***

Specify parameters $(\lambda,\theta)$.

**(1)** Generate $U$.

**(2)** If $U \le \frac{\theta}{\lambda+\theta}$, generate $X = G(1,\theta) = -\theta^{-1}\log(U)$; otherwise generate $X = G(2,\theta) = X_1 + X_2$, where $X_i = -\theta^{-1}\log(U_i), i = 1,2$. ($U_1 = U$ and $U_2$ is an additional uniform random variate.)

**(3)** Accept $Z = -\log(X)$ as a LL random variate (transform $X$ ).

**Remark 1**: As an analytic comparison, ***Algorithm 2*** requires two to three log function evaluations to generate one LL random variate while Jodra and Jimenez-Gamero [5] quantile approach needs one calculation of the Lambert W function.

***2.5 Survival and Hazard Rate functions***

$$S(x;\theta,\lambda,p)=\frac{(1+p+\theta\lambda)\{\Gamma(1+p)-(-\theta\log x)^{1+p}\,Ei(-p,-\theta\log x)\}-x^{\theta}(-\theta\log x)^{1+p}}{\Gamma(1+p)[1+p+\lambda\theta]}$$

$$=\frac{(1+p+\theta\lambda)\{\Gamma(1+p)-\Gamma(1+p,-\theta\log x)\}-x^{\theta}(-\theta\log x)^{1+p}}{\Gamma(1+p)[1+p+\lambda\theta]}$$

$$r(x;\theta,\lambda,p)=\frac{\theta^{2}(-\theta\log x)^{p}(\lambda-\log x)\,x^{\theta-1}}{(1+p+\theta\lambda)\{\Gamma(1+p)-(-\theta\log x)^{1+p}\,Ei(-p,-\theta\log x)\}-x^{\theta}\,(-\theta\log x)^{1+p}}$$

$$=\frac{\theta^{2}(\lambda-\log x)(-\theta\log x)^{p}}{(1+p+\theta\lambda)\{\Gamma(1+p)-\Gamma(1+p,-\theta\log x)\}-x^{\theta}\,(-\theta\log x)^{1+p}}$$

It can be easily checked that for $p=0$, the above results reduce to those of the $LL_0(\theta,\lambda)$ distribution of Gomez-Deniz et al. [3]. When $p$ is an integer, the expressions can be written in compact form using the result $\Gamma(p+1,y)=p!\,e^{-y}\sum_{k=0}^{p}y^{k}/k!$.

Some illustrative plots of the hazard function of $LL_p(\theta,\lambda)$ distributions for different choices of parameters $\theta,\lambda$ and index $p$ are presented in Figure 2, which reveals that the hazard function can be increasing and bath-tub shaped.

**Figure 2 here**

***2.6 Shannon Entropy of $LL_p(\theta,\lambda)$ and related results***

Entropy is a quantity that is often used in model selection. According to the maximum entropy principle, to make inference based on incomplete information, a distribution that best represents the current state of knowledge is the one that has the maximum entropy. Let $X$ be a random variable with pdf $f(x)$. The Shannon entropy of a pdf $f(x)$ is defined by

$$H(X) = H(f) = -\int f(x)\log f(x)\,dx .$$

We now derive the Shannon entropy of the GLL distribution, $LL_p(\theta,\lambda)$.

$$H_p(X) = E[-\log f(X;\theta,\lambda,p)]$$

$$= -E\left[\log\left[\frac{\theta^{2+p}}{\Gamma(1+p)[1+p+\lambda\theta]}(-\log X)^p(\lambda-\log X)X^{\theta-1}\right]\right]$$

$$= -\log\left[\frac{\theta^{2+p}}{\Gamma(1+p)[1+p+\lambda\theta]}\right] - pE\left[\log(\log(1/X)\right] - E[\log\{(\lambda-\log X)X^{\theta-1}\}] . \qquad (5)$$

It can be checked that for $p = 0$, $H_p(X)$ reduces to $H_0(X)$, where $H_0(X)$ is the Shannon entropy of the LL distribution, $LL_0(\theta,\lambda)$ given by

$$H_0(X) = \frac{1}{\theta(1+\theta\lambda)}\left[\begin{array}{l}\theta(1-\lambda)(1-\theta)+\theta e^{\lambda\theta}Ei(-\lambda\theta) \\ \qquad -\theta(1+\lambda\theta)\log\left(\dfrac{\lambda\theta^2}{1+\lambda\theta}\right)-2\end{array}\right]$$

(see proposition 2 in p. 51 of [3]). Here $Ei(z)$ [3,p.51] is defined as

$$Ei(z) = -\int_z^{\infty}\frac{e^{-\omega}}{\omega}\,d\omega .$$

For integral values of index parameter $p$, exact expression for $E\left[\log(\log(1/X))\right]$ and $E[\log\{(\lambda-\log X)X^{\theta-1}\}]$ can be obtained using Mathematica as follows.

***p = 1,***

$$E\left[\log(\log(1/X))\right] = \frac{3+\theta\lambda-(2+\theta\lambda)(\gamma+\log\theta)}{2+\theta\lambda},$$

$$E[\log\{(\lambda-\log X)X^{\theta-1}\}] = \frac{6-3\theta+2\theta\lambda(1-\theta)-e^{\theta\lambda}(2-\theta\lambda)Ei(-\theta\lambda)+\theta(2+\theta\lambda)\log\lambda}{\theta(2+\theta\lambda)}$$

$p = 2$,

$$E[\log\{\log(1/X)\}] = \frac{11 + 3\theta\lambda - 2(3+\theta\lambda)(\gamma + \log\theta)}{6 + 2\theta\lambda}$$

$$E[\log\{(\lambda - \log X)X^{\theta-1}\}] = \frac{24 - \theta(13 - 6\lambda + 7\theta\lambda) - e^{\theta\lambda}\theta(6 + \theta\lambda(\theta\lambda - 4))Ei(-\theta\lambda) + 2\theta(3+\theta\lambda)\log\lambda}{\theta(6 + 2\theta\lambda)}$$

and so on.

In particular for $p = 1$, Shannon Entropy for $LL_1(\theta, \lambda)$ using (5) is given by

$$H_1(X) = \log\left[\frac{(2+\lambda\theta)}{(1+\lambda\theta)\theta^2}\right] - \frac{3 + \theta\lambda - (2+\theta\lambda)(\gamma + \log\theta)}{2 + \theta\lambda}$$
$$- \frac{6 - 3\theta + 2\theta\lambda(1-\theta) - e^{\theta\lambda}(2 - \theta\lambda)Ei(-\theta\lambda) + \theta(2+\theta\lambda)\log\lambda}{\theta(2+\theta\lambda)}$$

where $\gamma$ is the Euler-Mascheroni constant.

Since the GLL distribution is a weighted LL distribution, it is of interest to examine the connection between weighted distributions and entropy. We first define weighted entropy (see, for example, [11]) and Kullback-Leibler distance.

**Definition 1**:

(a) For a pdf $f(x)$ and a nonnegative weight function $w(x)$, the weighted entropy of $f(x)$ is defined by

$$H(w; f) = -\int w(x) f(x) \log f(x) dx.$$

(b) The Kullback-Leibler distance between two pdf's $f(x)$ and $g(x)$ is defined by

$$KL(f; g) = \int f(x) \log\{f(x)/g(x)\} dx.$$

We state the connection between weighted distribution and entropy.

**Theorem 1**: Let $f(x)$ and $g(x)$ be two pdf's. The pdf $f(x)$ is a weighted $g(x)$ with weight $w(x)$, if and only if,

$$H(f) = -KL(f;g) + \frac{1}{E[w(X)]} H(w;g).$$

***Proof.*** If $f(x) = \frac{w(x)g(x)}{E[w(X)]}$, the proof follows directly from the definition of Shannon entropy.

Conversely, if the expression holds, straightforward manipulation leads to

$$\int f(x)\log g(x)dx - \frac{1}{E[w(X)]}\int w(x)g(x)\log g(x)dx = 0.$$

That is, pdf $f(x)$ is a weighted $g(x)$ with weight $w(x)$.

**Notes**.

**(1)** We may express Shannon entropy of $LL_p(\theta,\lambda)$ in terms of weighted entropy of $LL_{p-1}(\theta,\lambda)$ as follows:

$$H_p(X) = E[-\log f(X;\theta,\lambda,p)] = -E\left[\log\left[\frac{\theta(p+\lambda\theta)}{p(1+p+\lambda\theta)}(-\log X)f(X;\theta,\lambda,p-1)\right]\right]$$

$$= -\log\left[\frac{\theta(p+\lambda\theta)}{p(1+p+\lambda\theta)}\right] - E[\log\{-\log X\}] - E[\log[f(X;\theta,\lambda,p-1)]]$$

Since

$$f(X;\theta,\lambda,p) = \frac{\theta(p+\lambda\theta)}{p(1+p+\lambda\theta)}(-\log X)f(X;\theta,\lambda,p-1), \qquad (*)$$

$$E[\log[f(X;\theta,\lambda,p-1)]] = \int_0^1 \log[f(X;\theta,\lambda,p-1)]f(X;\theta,\lambda,p)dx$$

$$= \int_0^1 \log[f(X;\theta,\lambda,p-1)]\frac{\theta(p+\lambda\theta)}{p(1+p+\lambda\theta)}(-\log X)f(X;\theta,\lambda,p-1)dx$$

$$= \frac{\theta(p+\lambda\theta)}{p(1+p+\lambda\theta)}\int_0^1 (-\log X)\log[f(X;\theta,\lambda,p-1)]f(X;\theta,\lambda,p-1)dx$$

$$= -\log\left[\frac{\theta(p+\lambda\theta)}{p(1+p+\lambda\theta)}\right] - E\left[\log\{-\log X\}\right]$$

$$-\frac{\theta(p+\lambda\theta)}{p(1+p+\lambda\theta)}\int_0^1 (-\log x)\log\left[f(X;\theta,\lambda,p-1)\right] f(X;\theta,\lambda,p-1)\,dx$$

$$= -\log\left[\frac{\theta(p+\lambda\theta)}{p(1+p+\lambda\theta)}\right] - E\left[\log\{-\log X\}\right] - \frac{\theta(p+\lambda\theta)}{p(1+p+\lambda\theta)} H\left(-\log x; f(X;\theta,\lambda,p-1)\right),$$

where $H\left(-\log x; f(X;\theta,\lambda,p-1)\right)$ is the weighted entropy of $LL_{p-2}(\theta,\lambda)$ with weight function ($-\log x$).

**(2)** By using (*) we obtain interesting expressions for moments $E\left[(\log X)^r\right]$ of $LL_p(\theta,\lambda)$.

Replace $p$ by $p+1$ in (*) and integrate both sides. We get

$$E(\log X;\theta,\lambda,p) = -\frac{(p+1)(2+p+\lambda\theta)}{\theta(1+p+\lambda\theta)}.$$

In general, $E\left[(\log X)^r;\theta,\lambda,p)\right] = -\dfrac{(p+1)_{[r]}(1+r+p+\lambda\theta)}{\theta^r(1+p+\lambda\theta)}.$

For $p=0$ gives the corresponding expression for LL distribution [3,equation (5)] as

$$E(\log X) = -\frac{(2+\lambda\theta)}{\theta(1+\lambda\theta)}.$$

**Remark 2**. Note that equation (5) of [3] contains a typo.

For $n$ $LL_p(\theta,\lambda)$ random variables $X_1$, $X_2$, …, $X_n$, we get

$$E(\log\prod_{i=1}^{n} X_i;\theta,\lambda,p) = -\frac{n(p+1)(2+p+\lambda\theta)}{\theta(1+p+\lambda\theta)} = -\frac{(p+1)\left(2+p+n\lambda\dfrac{\theta}{n}\right)}{\dfrac{\theta}{n}\left(1+p+n\lambda\dfrac{\theta}{n}\right)}$$

which is $E(\log Y;\theta,\lambda,p)$, where $Y \sim LL_p\left(\dfrac{\theta}{n}, n\lambda, p\right)$. Additionally, we can get

$$E(x\log X;\theta,\lambda,p)=-\left(\frac{\theta}{1+\theta}\right)^{p+3}\frac{(p+1)(2+p+\lambda(1+\theta))}{\theta(1+p+\lambda\theta)},\ldots,$$

$$E(X^{r}\log X;\theta,\lambda,p)=-\left(\frac{\theta}{r+\theta}\right)^{p+3}\frac{(p+1)(2+p+\lambda(r+\theta))}{\theta(1+p+\lambda\theta)}.$$

***2.7 Stochastic Ordering***

The comparison of random quantities through the notion of stochastic ordering has important applications in many areas, for instance, in risk theory and reliability theory, see Chapters 12, 15 and 16 in [12]. We will consider the likelihood ratio (*LR*), hazard rate (*HR*) and stochastic (*ST*) orderings for $LL_p(\theta,\lambda)$ random variables in this section. The definition of likelihood ratio, stochastic and hazard rate orders are as follow.

**Definition 2**: Suppose two random variables *X* and *Y* have pdf's *f* and *g*, cdf's $F$ and $G$, hazard rate $h_X$ and $h_Y$ respectively. Then

**(1)** *X* is said to be smaller than *Y* in the ***likelihood ratio order***, denoted by $X\le_{LR}Y$, if $f(x)g(y)\ge f(y)g(x)$ for all $x\le y$.

**(2)** *X* is said to be ***stochastically*** smaller than *Y*, denoted by $X\le_{ST}Y$, if $F(x)\ge G(x),\forall x$.

**(3)** *X* is said to be smaller than *Y* in ***hazard rate order***, denoted by $X\le_{HR}Y$, if $h_X(x)\le h_Y(x)$ for all $x$.

**Theorem 2**. Let $X_1$ and $X_2$ be random variables following $LL_{p_1}(\theta_1,\lambda_1)$ and $LL_{p_2}(\theta_2,\lambda_2)$ distributions, respectively. If $\theta_1\le\theta_2,\lambda_1\le\lambda_2$ and $p_2\le p_1$ then $X_1\le_{LR}X_2$.

***Proof***: Consider the ratio

$$\frac{f(x;\theta_2,\lambda_2,p_2)}{f(x;\theta_1,\lambda_1,p_1)}=\frac{\theta_2^{\,2+p_2}\Gamma(1+p_1)[1+p_1+\lambda_1\theta_1]}{\theta_1^{\,2+p_1}\Gamma(1+p_2)[1+p_2+\lambda_2\theta_2]}(-\log x)^{p_2-p_1}h(x)\ ,\qquad (6)$$

where $h(x)=\dfrac{\lambda_2-\log x}{\lambda_1-\log x}x^{\theta_2-\theta_1}$. Gomez-Deniz et al. [3] have shown that the function $h(x)$ is non-decreasing for $x\in(0,1)$ if $\theta_1\le\theta_2,\lambda_1\le\lambda_2$. If $p_2\le p_1$, it is clear that $(-\log x)^{p_2-p_1}$ is non-decreasing for $x\in(0,1)$. This implies that if $\theta_1\le\theta_2,\lambda_1\le\lambda_2$ and $p_2\le p_1$, then the ratio in (6) is non-decreasing for $x\in(0,1)$ and hence, $X_1\le_{LR}X_2$. □

Clearly, Theorem 2 may be applied in the context of monotone likelihood ratio to obtain uniformly most powerful tests for parameters of interest.

LR ordering is stronger than hazard rate and stochastic orderings and this leads to the following implications [3]: $X_1\le_{LR}X_2\Rightarrow X_1\le_{HR}X_2\Rightarrow X_1\le_{ST}X_2$. Therefore, as in Corollary 1 of Gomez-Deniz et al. [3] similar results can be shown for the GLL distribution as follows.

**Corollary 1**. Let $X_1$ and $X_2$ be random variables following $LL_{p_1}(\theta_1,\lambda_1)$ and $LL_{p_2}(\theta_2,\lambda_2)$ distributions, respectively. If $\theta_1\le\theta_2,\lambda_1\le\lambda_2$ and $p_2\le p_1$, then

**(1)** the moments, $E[X_1^k]\le E[X_2^k]$ for all $k>0$;

**(2)** the hazard rates, $r_1(x)\le r_2(x)$ for all $x\in(0,1)$.

**Remark 3.** For a nonnegative random variable $X$, let $X^w$ be the corresponding weighted random variable derived using a nonnegative weight function $w(x)$ (see section 2), then $X\underset{lr}{\ge}X^w$ if $w(x)$ is decreasing [13]. In view of the results in section 2, that $LL_p(\theta,\lambda)$ is a weighted $LL\ (\theta,\lambda)$ with nonnegative decreasing weight function $w(x)=(-\log x)^p$ it is therefore obvious that $LL(\theta,\lambda)\ge_{LR}LL_p(\theta,\lambda),\ p>0$. In fact, it can be further generalized to state that $LL_p(\theta,\lambda)\ge_{LR}LL_q(\theta,\lambda),\ q>p$.

The corresponding hazard rate and stochastic orderings follow as a consequence.

*2.8 Convexity, Concavity and Log-concavity for GLL Distribution*

Gomez-Deniz et al. [3] have shown the LL cdf $F(x)$ to be convex but did not consider the property of log-concavity. In this section we examine the convexity of the cdf and log-concavity of the pdf for the GLL distribution. For brevity, we suppress the cdf notation for $LL_p(\theta,\lambda)$ distribution in (3) as $F(x)$ in this section.

**Theorem 3**. If $0<\theta\le 1$, $p\ge 0$ and $\lambda\ge 0$, then $F(x)$ is concave for $x\in(0,1)$. Hence, for $0<\theta\le 1$, $F(x)$ is also log-concave for $x\in(0,1)$ since $F(x)\ge 0$.

***Proof***: If $0<\theta\le 1$, then $(-\log x)^p$, $(\lambda-\log x)$ and $(x)^{\theta-1}$ are decreasing in $x\in(0,1)$ for $p\ge 0$ and $\lambda\ge 0$. This implies that the pdf

$$F'(x)=f(x;\theta,\lambda,p)=\frac{\theta^{2+p}}{\Gamma(1+p)[1+p+\lambda\theta]}(-\log x)^p(\lambda-\log x)\,x^{\theta-1},$$

is decreasing in $x\in(0,1)$. Thus, $F(x)$ is concave for $x\in(0,1)$.

Since $F(x)\ge 0$, concavity implies $F(x)$ is also log-concave for $x\in(0,1)$. □

**Theorem 4**. The function $F(x)$ is neither convex nor concave for $x\in(0,1)$ for any $\theta>1$, $p>0$ and $\lambda\ge 0$.

***Proof***: For any $\theta>1$, $p>0$ and $\lambda\ge 0$, consider the second order derivative of $F(x)$ given by

$$F''(x)=\frac{\theta^{2+p}}{\Gamma(1+p)[1+p+\lambda\theta]}x^{\theta-1}(-\log x)^{p-1}(\theta-1)\left\{(\log x)^2-\left(\lambda-\frac{p+1}{\theta-1}\right)\log x-\frac{p\lambda}{\theta-1}\right\}.$$

In order for $F(x)$ to be convex, $F''(x)$ must be $\ge 0$ or that the term $\left\{(\log x)^2-\left(\lambda-\frac{p+1}{\theta-1}\right)\log x-\frac{p\lambda}{\theta-1}\right\}\ge 0$ for all $x\in(0,1)$. However, when $x\to 1$,

$$\left\{(\log x)^2-\left(\lambda-\frac{p+1}{\theta-1}\right)\log x-\frac{p\lambda}{\theta-1}\right\}\to-\frac{p\lambda}{\theta-1}\le 0.$$

Thus, $F(x)$ is not convex for $x\in(0,1)$ for any $\theta>1$, $p>0$ and $\lambda\ge 0$.

When $x\to 0$,

$$\left\{(\log x)^2-\left(\lambda-\frac{p+1}{\theta-1}\right)\log x-\frac{p\lambda}{\theta-1}\right\}=(\log x)^2\left\{1-\frac{\left(\lambda-\frac{p+1}{\theta-1}\right)}{\log x}-\frac{p\lambda}{(\theta-1)(\log x)^2}\right\}\to\infty.$$

This implies that $F''(x)$ cannot be $\le 0$ for all $x\in(0,1)$. Therefore, $F(x)$ is not concave for $x\in(0,1)$ for any $\theta>1$, $p>0$ and $\lambda\ge 0$. $\square$

**Note**: $F(x)$ is convex for $x\in(0,1)$ only for $p=0$ [that is, when $LL_p(\theta,\lambda)$ distribution reduces to the LL distribution] and $\lambda(\theta-1)\ge 1$ as shown in Theorem 4 of [3].

We next show that the $LL_p(\theta,\lambda)$ distribution is log-concave. A function $f(x)$ is log-concave if $\log f(x)$ is concave and log-convex if $\log f(x)$ is a convex function. Based on Definition 2 of Borzadaran and Borzadaran [14], the log-concavity of a function $f(x)$ on an interval $(a,b)$ is equivalent to $f'(x)/f(x)$ being monotonically decreasing in $(a,b)$ or $\left(\ln f(x)\right)''<0$ .

**Theorem 5**. If $\theta>1$, the $LL_p(\theta,\lambda)$ pdf is log-concave.

***Proof.*** It is sufficient to show that

$$\frac{f'(x)}{f(x)}=\frac{F''(x)}{F'(x)}=\frac{(\theta-1)}{(-\log x)}\left\{(\log x)^2-\left(\lambda-\frac{p+1}{\theta-1}\right)\log x-\frac{p\lambda}{\theta-1}\right\}$$

is monotonically decreasing in $x\in(0,1)$. Rewrite the above as

$$\frac{f'(x)}{f(x)} = (\theta-1)\left\{(-\log x) - \left(\lambda - \frac{p+1}{\theta-1}\right) - \frac{p\lambda}{\theta-1}\frac{1}{(-\log x)}\right\}.$$

Noting that $\log x$ is monotonically increasing in $(0,1)$, $-\log x$ is monotonically decreasing. For $\theta > 1$, $f'(x)/f(x)$ is monotonically decreasing in $(0,1)$. Thus $f(x)$ is log-concave. □

**Remark 4**. We have also shown that the LL ($LL_0(\theta,\lambda)$) pdf is log-concave provided $\theta > 1$.

The log-concavity or log-convexity of pdfs implies many interesting properties of the distributions, especially reliability properties. There are many applications of this log-concavity property in diverse disciplines [15]. Many properties of the GLL distribution follow from the property of log-concavity [14,p.205–206 ]. For instance,

**(1)** GLL pdf is strongly unimodal; a distribution $F$ on $R$ is said to be strongly unimodal if the convolution of $F$ with any unimodal distribution is again unimodal.

**(2)** GLL pdf is a Polya frequency density of order 2, that is,

$f(x-y)f(x'-y') - f(x-y')f(x'-y) \geq 0$ for $x < x', y < y'$.

**(3)** Hazard function is a non-decreasing function in $x$.

**(4)** Distribution function and survival function are log-concave.

**3. Application to insurance premium loading**

In this section we apply the results of Section 2.8 to insurance premium loading. According to the basic premium principle, assuming common agreement on the risk distribution, the net premium $P(X)$ to be charged on an insurance coverage for exposure to risk $X$ is given by $P(X) = E[X]$. If there is no agreement on the risk distribution, a loading is added to $X$. One approach to add this loading is through transformation of the initial cumulative distribution function of $X$ by a continuous and non-decreasing function $h$ known as the distortion function.

This transformation results in a new distribution corresponding to a random variable $Y$. If the distortion function $h$ is convex, this guarantees that $X \leq_{ST} Y$ and which further implies that $E[X] \leq E[Y]$, that is, loading is nonnegative (refer to [3] for further discussion).

**Theorem 6.** If $F(x;\theta,\lambda,p)$ of $LL_p(\theta,\lambda)$ given by (3) is concave, then $1-F(1-x;\theta,\lambda,p)$ is a convex function from (0, 1) to (0, 1) for $0<\theta\leq 1$ and $p\geq 0, \lambda\geq 0$.

***Proof***: It has been shown in Theorem 3 that the cdf $F(x;\theta,\lambda,p)$ of $LL_p(\theta,\lambda)$ is concave for $0<\theta\leq 1$ and $p\geq 0,\lambda\geq 0$. Hence, $1-F(1-x;\theta,\lambda,p)$ is a convex function from (0, 1) to (0, 1).

**Remark 5**: $F(x;\theta,\lambda,p)$ can be used as a distortion function to distort survival function (sf) of a given random variable as stated in the Corollary 2 next.

**Corollary 2.** If $X$ is the risk with sf $\overline{G}(x)$ and let $Z$ be a distorted random variable with sf $F[\overline{G}(x);\theta,\lambda,p]$ for $0<\theta<1$ and $p\geq 0,\lambda\geq 0$. Let $E[Z]=\int_0^\infty F[\overline{G}(x);\theta,\lambda,p]dx = P_{\theta,\lambda}(X)$ be the distorted premium and $P_n(X)=\int_0^\infty [\overline{G}(x)]^n\, dx$, $0<n\leq 1$ be the proportional hazard premium (see [7]) of an insurance product respectively. Then $P_{\theta,\lambda}(X)$ is a premium principle such that

**(1)** $P_n(X)\leq P_{\theta,\lambda}(X)\leq \max(X)$, for all $n\geq\theta$,

**(2)** $P_{\theta,\lambda}(aX+b)=aP_{\theta,\lambda}(X)+b$,

**(3)** if $\overline{G}_1(X_1)$ and $\overline{G}_1(X_2)$ are sf of two non-negative risk random variables $X_1$ and $X_2$ with $\overline{G}_1(X_1)\leq\overline{G}_2(X_2)$, that is, $X_1$ precedes $X_2$ under first stochastic dominance then $P_{\theta,\lambda}(X_1)\leq P_{\theta,\lambda}(X_2)$,

**(4)** if $X_1$ precedes $X_2$ under second stochastic dominance, that is, if

$\int_x^{\infty} \overline{G}_1(x_1)\,dx_1 \le \int_x^{\infty} \overline{G}_2(x_2)\,dx_2$ for all $x \ge 0$ then $P_{\theta,\lambda}(X_1) \le P_{\theta,\lambda}(X_2)$.

***Proof***: From Theorem 3 above we know that $F(x;\theta,\lambda,p)$ is concave for $x \in (0,1)$ when $0 < \theta \le 1$, $p \ge 0$ and $\lambda \ge 0$. Also being a cdf, it is an increasing function of $x$ with $F(0;\theta,\lambda,p) = 0$ and $F(1;\theta,\lambda,p) = 1$.

The results (2), (3) and (4) follow immediately from Definition 6 of distortion premium principle and subsequent properties thereof in [8]. For the result (1) from [8] it follows that $P_{\theta,\lambda}(X) \le \max(X)$.

Now we provide a proof of $P_n(X) \le P_{\theta,\lambda}(X)$. We first prove that for any $p$, $F(x;\theta,\lambda,p) \ge x^{\theta}$.

Case I: When $p$ is an integer, we get

$$F(x;\theta,\lambda,p) = \frac{x^{\theta}\theta^{1+p}(-\log x)^{1+p} + (1+p+\theta\lambda)\,p!\,x^{\theta}\sum_{k=0}^{p}(-\theta\log x)^k/k!}{(1+p+\lambda\theta)\,p!}$$

$$= x^{\theta}\left[1 + \sum_{k=1}^{p}\frac{(-\theta\log x)^k}{k!} + \frac{(-\theta\log x)^{1+p}}{(1+p+\lambda\theta)\,p!}\right] > x^{\theta}.$$

Case II: For real $p$, we apply the integral representation (https://en.wikipedia.org/wiki/ Incomplete_gamma_function; see section "Evaluation Formulae")

$$\Gamma(1+p,y) = e^{-y}\,y^{1+p}\int_0^{\infty} e^{-yu}\left(1+u\right)^p du$$

$$> e^{-y}\,y^{1+p}\int_0^{\infty} e^{-yu}\left(u\right)^p du = e^{-y}\int_0^{\infty} e^{-t}t^{(p+1)-1}dt = e^{-y}\Gamma(1+p)$$

Hence,

$$F(x;\theta,\lambda,p)=\frac{x^{\theta}(-\theta\log x)^{1+p}+(1+p+\theta\lambda)\Gamma(1+p,-\theta\log x)}{\Gamma(1+p)[1+p+\lambda\theta]}$$

$$>\frac{x^{\theta}(-\theta\log x)^{1+p}+(1+p+\theta\lambda)e^{\theta\log x}\Gamma(1+p)}{\Gamma(1+p)[1+p+\lambda\theta]}=x^{\theta}\left[1+\frac{(-\theta\log x)^{1+p}}{\Gamma(1+p)[1+p+\lambda\theta]}\right]>x^{\theta}$$

Therefore, for any $p$, $F(x;\theta,\lambda,p) \geq x^{\theta}$. It follows that

$$F(x;\theta,\lambda,p)\geq x^{\theta} \Rightarrow x^{\theta}\leq F(x;\theta,\lambda,p) \text{ for } x\in(0,1)$$

$$\Rightarrow (\overline{G}(x))^{\theta}\leq F(\overline{G}(x);\theta,\lambda,p) \text{ since } 0\leq\overline{G}(x)\leq 1 \text{ for all } x\in\Re.$$

$$\Rightarrow (\overline{G}(x))^{\theta}\leq F(\overline{G}(x);\theta,\lambda,p) \text{ for all } x\in\Re.$$

Now, for all $x\in\Re$, we have

$$(\overline{G}(x))^{n}\leq F(\overline{G}(x);\theta,\lambda,p) \text{ when } 0<\theta\leq n<1$$

$$\Rightarrow \int_0^{\infty}(\overline{G}(x))^{n}dx\leq\int_0^{\infty}F(\overline{G}(x);\theta,\lambda,p))dx \Rightarrow P_n(X)\leq P_{\theta,\lambda}(X).$$

**Remark 6.** This new distorted premium principle is a trade-off between the proportional hazard premium and maximal premium. For $n\to 1$, $P_n(X)\to E(X)$ and for $n\to 0$, $P_n(X)\to \max(X)$ (see [7]).

For a numerical confirmation of the Corollary 2 (1) we present in Table 1, values of the proportional hazard premium $P_n(X)$ and the distorted premium $P_{\theta,\lambda}(X)$ obtained using the Generalized Log–Lindley distribution for different values of the parameters considering the Exponential, Weibull, and Inverse Gaussian as the underlying risk distribution.

**Table 1 here**

## 4. A useful re-parameterization of $LL_p(\theta,\lambda)$

Starting with a two-parameter Lindley distribution [16], Jodra et al. [5] obtained a re-parameterized version of (1) with pdf

$$f(x;\theta,\pi)=\theta(\pi+\theta(\pi-1)\log x)\,x^{\theta-1}, 0<x<1,\ 0\leq\pi\leq 1,\ \theta>0. \qquad (7)$$

In fact, this can be obtained by substituting $\frac{\lambda\theta}{1+\lambda\theta}=\pi$ that is $\frac{1}{1+\lambda\theta}=1-\pi$ in the LL distribution in (1).

The pdf in (7) overcomes the issue of unbounded parameter space of (1). Applying the same re-parameterization to $LL_p(\theta,\lambda)$ distribution in (2) we obtain a re-parameterized version with pdf

$$f(x;\theta,\pi,p)=\frac{\theta^{1+p}}{\Gamma(1+p)[1+(1-\pi)p]}(-\log x)^p(\pi+\theta(\pi-1)\log x)\,x^{\theta-1}, \qquad (8)$$

where $,0<x<1,\ 0\le\pi\le1,\ \theta>0, p\ge0$.

Now the mean of the re-parameterized GLL can be easily derived using (4) as

$$E(X;\theta,\pi,p)=\left(\frac{\theta}{1+\theta}\right)^{2+p}\frac{\pi+\theta+p\theta(1-\pi)}{\theta+p\theta(1-\pi)}=\left(\frac{\theta}{1+\theta}\right)^{2+p}\frac{1+p(1-\pi)+\pi/\theta}{1+p(1-\pi)}.$$

For $p$ = 0, we get back (1) from (7) and (2) from (8). Note that the new parameter $\pi$ introduced in (8) is bounded.

**5. Maximum Likelihood Estimation and Information Matrix**

The likelihood function for a random sample of size $n$ from the $LL_p(\theta,\lambda)$ is

$$L=\frac{\theta^{n(2+p)}}{\{\Gamma(1+p)(1+p+\lambda\theta)\}^n}\left(\prod_{i=1}^{n}-\log x_i\right)^p\prod_{i=1}^{n}(\lambda-\log x_i)\,x^{\theta-1}.$$

The log-likelihood function is then given by

$$l=\log L=n(2+p)\log\theta-n\log\Gamma(1+p)-n\log(1+p+\lambda\theta)+p\sum_{i=1}^{n}\log(-\log x_i)$$

$$+\sum_{i=1}^{n}\log(\lambda-\log x_i)+(\theta-1)\sum_{i=1}^{n}\log x_i.$$

The first and second order derivatives of the log-likelihood function are:

$$\frac{\partial l}{\partial \theta}=\frac{n(2+p)}{\theta}-\frac{n\lambda}{1+p+\lambda\theta}+\sum_{i=1}^{n}\log x_i$$

$$\frac{\partial l}{\partial \lambda}=\frac{n\theta}{1+p+\lambda\theta}+\sum_{i=1}^{n}\frac{1}{\lambda-\log x_i}$$

$$\frac{\partial l}{\partial p}=n\log\theta-n\frac{\Gamma^{/}(1+p)}{\Gamma(1+p)}-\frac{n}{1+p+\lambda\theta}+\sum_{i=1}^{n}\log[\log[1/x_i]]$$

$$\frac{\partial^2 l}{\partial \theta^2}=\frac{-n(2+p)}{\theta^2}+\frac{n\lambda^2}{(1+p+\lambda\theta)^2}=-\frac{n[(2+p)(1+p+\lambda\theta)^2-\lambda^2\theta^2]}{\theta^2(1+p+\lambda\theta)^2}$$

$$\frac{\partial^2 l}{\partial \lambda^2}=-\frac{n\theta^2}{(1+p+\lambda\theta)^2}-\sum_{i=1}^{n}\frac{1}{(\lambda-\log x_i)^2}$$

$$\frac{\partial^2 l}{\partial p^2}=n\left[\frac{\Gamma^{/}(1+p)}{\Gamma(1+p)}\right]^2-n\frac{\Gamma^{//}(1+p)}{\Gamma(1+p)}$$

$$\frac{\partial^2 l}{\partial p\,\partial \theta}=\frac{n}{\theta}+\frac{n\lambda}{(1+p+\lambda\theta)^2}$$

$$\frac{\partial^2 l}{\partial \lambda\,\partial \theta}=\frac{n\lambda\theta}{(1+p+\lambda\theta)^2}-\frac{n}{1+p+\lambda\theta}$$

$$\frac{\partial^2 l}{\partial p\,\partial \lambda}=-\frac{n\theta}{(1+p+\lambda\theta)^2}$$

For information matrix, we obtain the following result

$$E\left[\sum_{i=1}^{n}\frac{1}{(\lambda-\log x_i)^2}\right]=\frac{\theta^{2+p}e^{\lambda\theta}}{\Gamma(1+p)[1+p+\lambda\theta]}\sum_{r=0}^{p}\binom{p}{r}\left(\frac{1}{\lambda\theta}\right)^{r}(-\lambda)^{p}\Gamma(r,\lambda\theta)\,.$$

For $p=0$, this reduces to the result given in [3].

The information matrix is given by $(I_{j,k}),\ j,k=1,2,3.$ Where

$$I_{11}=\frac{n\{(2+p)(1+p+\lambda\theta)^2-\lambda^2\theta^2\}}{\theta^2(1+p+\lambda\theta)^2},$$

$$I_{12} = -\frac{n\lambda\theta}{(1+p+\lambda\theta)^2} + \frac{n}{1+p+\lambda\theta}, I_{13} = -\frac{n}{\theta} - \frac{n\lambda}{(1+p+\lambda\theta)^2},$$

$$I_{22} = -\frac{n\theta^2}{(1+p+\lambda\theta)^2} + \frac{\theta^{2+p} e^{\lambda\theta}}{\Gamma(1+p)[1+p+\lambda\theta]} \sum_{r=0}^{p} \binom{p}{r} \left(\frac{1}{\lambda\theta}\right)^r (-\lambda)^p \Gamma(r,\lambda\theta)$$

$$I_{23} = \frac{n\theta^2}{(1+p+\lambda\theta)^2} \text{ and } I_{33} = -n\left\{\frac{\Gamma^{/}(1+p)}{\Gamma(1+p)}\right\}^2 + n\frac{\Gamma^{//}(1+p)}{\Gamma(1+p)}.$$

This matrix can be inverted to get the asymptotic variance-covariance matrix for the maximum likelihood estimates.

The derivatives of the gamma function $\Gamma(\alpha)$ in $I_{33}$ may be computed as follows. Let $\psi(\alpha) = \frac{d \log\Gamma(\alpha)}{d\alpha}$ and $\psi'(\alpha) = \frac{d^2 \log\Gamma(\alpha)}{d\alpha^2}$ be the digamma and trigamma functions, respectively. We note the following formulae for the derivatives of the gamma function $\Gamma(\alpha)$ in terms of digamma and trigamma functions:

$$\Gamma'(\alpha) = \Gamma(\alpha)\psi(\alpha),\ \Gamma''(\alpha) = \Gamma(\alpha)\left\{\psi'(\alpha) + \psi^2(\alpha)\right\}$$

Spouge [17] has given efficient methods to compute the digamma and trigamma functions. In the programming language *R* the functions *digamma(x)* and *trigamma(x)* are available to compute these functions.

**6. Data Modeling Applications**

In this section, parameter estimation was performed using the R packages. Beta regression is performed using the package *betareg* with logit link for the mean, while other regressions are performed through optimization using the function mle2 in the package *bbmle*.

***6.1 Risk Management Data***

In this section, we consider the modeling of the data set originally used in [6] and considered by Gomez-Deniz et al. [3] about the cost effectiveness of risk management (measured in percentages) in relation to exposure to certain property and casualty losses, adjusted by several other variables such as size of assets and industry risk. Description of the data set may be found in [3]. We take the response variable to be $Y$ = FIRMCOST/100. Six other variables (covariates) are ASSUME ($X_1$), CAP ($X_2$), SIZELOG ($X_3$), INDCOST ($X_4$), CENTRAL ($X_5$) and SOPH ($X_6$). We model response variable ($Y$) as well as its complimentary (1-$Y$) without and with covariates. To accommodate the covariates we introduce two regression models by first linking the covariates to the parameter $\theta$ through the log link function, and then to the mean $\mu$ through the logit link function for both the LL and GLL distributions. It may be noted that Jodra and Jimenez-Gamero [5] also investigated the same data set using a re-parameterization of the LL distribution.

*6.1.1 Modeling the parameter centering parameter $\theta$ of $LL_p(\theta,\lambda)$*

Here we first note that as $\lambda \to \infty$, the mean of the $LL_p(\theta,\lambda)$ distribution in (4),

$$E(X;\theta,\lambda,p)=\left(\frac{\theta}{1+\theta}\right)^{2+p}\frac{1+p+\lambda(1+\theta)}{1+p+\lambda\theta}\to\left(\frac{\theta}{1+\theta}\right)^{1+p}.$$

This implies that the parameter $\theta$ plays a certain "centering" role for the distribution. For the purpose of regression modeling of the parameter $\theta$, a suitable link function is required. Suppose that a random sample $Y_1,Y_2,\ldots,Y_n$ of size $n$ is obtained from the $LL_p(\theta,\lambda)$ distribution. For a set of $k$ covariates, the log link for the $LL_p(\theta,\lambda)$ regression model gives the parameter $\theta$ for each $Y_i$ as

$$\theta_i=\frac{\exp(\mathbf{x}_i^T\boldsymbol{\beta})}{1+\exp(\mathbf{x}_i^T\boldsymbol{\beta})},\quad i=1,2,\ldots,n,$$

where $\mathbf{x}_i^T=(1,x_{i1},\ldots,x_{ik})$ are the covariates with corresponding coefficients $\boldsymbol{\beta}=(\beta_0,\beta_1,\ldots,\beta_k)$.

Here, the beta (with logit link to mean), LL (with log link to $\theta$) and $LL_p(\theta,\lambda)$ (with log link to $\theta$) regression models are considered and the log-likelihood values and parameter estimates for the models considered, without and with covariates, are presented in Table 2.

**Table 2 here**

In terms of the log-likelihood values, it is clear from the results in Table 2 that the generalized Log-Lindley model fits the data best with or without covariates for the response variable *Y*. It is also the best model for the regression of 1 – *Y*, while the beta model is the best for this case without covariates. For the case of modeling 1 – *Y* with and without covariates, it is seen that the estimates for the parameter *p* of the $LL_p(\theta,\lambda)$ model approaches 0 and hence, approaches the results for the LL distribution.

The estimates for Log-Lindley regression model may not be compared directly with that in [3] since different parameterization is used. The estimates for the Log-Lindley model with response 1 – *Y* obtained here is close to that obtained in [5]. For both responses *Y* and 1 – *Y*, it is clear that the covariates SIZELOG and INDCOST are statistically significant across level of significance > 0.02. Hence, the measure of the firm's risk management cost effectiveness (FIRMCOST) is negatively associated with the logarithm of total assets (SIZELOG) but positively associated with the measure of the firm's industry risk (INDCOST).

*6.1.2 Modeling in terms of Mean*

Here we have investigated the cases of *Y* and 1 – *Y* without and with covariates for $LL_p(\theta,\lambda)$ re-parameterized in terms of its mean, $\mu$ in (4) together with two other new parameters $\phi$ and $\gamma$ such that the parameters in (2) are replaced by:

$$\theta=\frac{\mu\gamma(2+\phi)+\sqrt{\mu^2\gamma^2\phi^2+4\mu\gamma(1+\phi)}}{2(1-\mu\gamma)(1+\phi)},\ p=\frac{\log\gamma}{\log[(1+\theta)/\theta]},\ \lambda=\left(\frac{1+p}{1+\theta}\right)\phi$$

The new parameters are such that $0 < \mu < 1$, $\phi > 0$, $1 \leq \gamma < 1/\mu$. The re-parametrized distribution will now be denoted in terms of the new parameters as $LL_{\gamma}(\mu,\phi)$. Clearly, when $\gamma$ = 1, the re-parameterized log-Lindley distribution of Gomez-Deniz et al. [3] is obtained. Now for a random sample $Y_1, Y_2, \ldots, Y_n$ of size $n$ from the $LL_{\gamma}(\mu,\phi)$ distribution and a set of $k$ covariates, the logit link for the $LL_{\gamma}(\mu,\phi)$ regression model gives the mean for each $Y_i$ as

$$\mu_i = \frac{\exp(\mathbf{x}_i^T \boldsymbol{\beta})}{1+\exp(\mathbf{x}_i^T \boldsymbol{\beta})}, \quad i = 1,2,\cdots,n.$$

Here, only LL and $LL_{\gamma}(\mu,\phi)$ regression models with the logit link to mean as above are applied to the same data set used. The log-likelihood values and parameter estimates for the models considered, with and without covariates, are presented in Table 3.

**Table 3 here**

It is clear that in terms of the log-likelihood values in Table 3 that $LL_{\gamma}(\mu,\phi)$ model fits the data better for the response variable $Y$ as well as $1 - Y$, for the regression of $1 - Y$, the estimates of the parameter $\gamma$ for the $LL_{\gamma}(\mu,\phi)$ model approaches 1 which is why results approach that of the LL regression model as depicted.

### *6.2 Outpatient Health Expenditures Data*

This data set is obtained with description of the variables from Frees [18]. We considered the outpatient health expenditures for 1,352 individuals in '000,000 USD, that is the response variable $Y$ = EXPENDOP/1,000,000. Five other variables (covariates) considered are AGE, GENDER, insure, RACE and PHSTAT. We model the response variable ($Y$) without and with covariates. Similarly, the covariates are linked to the parameter $\theta$ for the log-Lindley and generalized log-Lindley regression models using the log link function in Section 6.1.1. The numerical results are presented in Table 4.

**Table 4 here**

Note that, from Table 4, the estimated parameter $\lambda$ for the log-Lindley distribution approached 0, which indicates that the model may be reduced to the special case corresponding with the transformed gamma distribution as noted in [3]. It is clear that in terms of the log-likelihood values in Table 4 that $LL_p(\theta, \lambda)$ model fits the data better for the response variable *Y* for both the case without and with covariates. From the estimates of all regression models, the variables AGE, GENDER, insure and PHSTAT are statistically significant across level of significance > 0.02. These variables are all positively associated with the amount of outpatient health expenditures *Y*.

## 7. Concluding discussion

A new distribution in (0, 1) is proposed that nests the Log-Lindley distribution of Gomez-Deniz et al. [3] and offers compact expressions for cdf and moments. This new distribution is shown to be a weighted Log-Lindley distribution which enhances its usefulness in statistical modeling. Many of its important structural properties like log-concavity and stochastic ordering are studied. An interesting characterization of the weighted distribution in terms of Kullback-Liebler distance and weighted entropy has been given. This is of utility for statistical inference with the proposed distribution as a weighted Log-Lindley distribution. A new class of distorted premium principle based on the proposed distribution is also introduced stating some important results. A re-parameterization of the proposed distribution is also prescribed for achieving bounded range for a parameter that allows regression of this parameter on covariates through a Logit link function. Application to two real life data sets, with a much better fit than the beta and Log-Lindley distributions, shows the relevance of the newly proposed distribution in modeling without covariates and also in regression analysis to accommodate covariates.

**Acknowledgements:**

**Table 1**: Values of proportional hazard premium, $P_n(X)$ and distorted premium, $P_{\theta,\lambda}(X)$ for fixed $p$ and for different distributions with given parameters

| | $P_n(X)$ | | | $P_{\theta,\lambda}(X)$, $p=1$ | | | | $P_{\theta,\lambda}(X)$, $p=2$ | | | |
|---|---|---|---|---|---|---|---|---|---|---|---|
| | $n$ | | | $(\theta,\lambda)$ | | | | $(\theta,\lambda)$ | | | |
| | 0.4 | 0.75 | 1.0 | (0.3, 0.5) | (0.3, 1.5) | (0.7, 0.5) | (0.7, 1.5) | (0.3, 0.5) | (0.3, 1.5) | (0.7, 0.5) | (0.7, 1.5) |
| **Exponential ($\lambda$ = rate)** | | | | | | | | | | | |
| $\lambda=0.5$ | 5.00 | 2.67 | 2.00 | 19.54 | 18.78 | 8.15 | 7.59 | 26.35 | 25.80 | 11.13 | 10.69 |
| $\lambda=2.0$ | 1.25 | 0.67 | 0.50 | 4.88 | 4.69 | 2.04 | 1.90 | 6.59 | 6.45 | 2.78 | 2.67 |
| **Weibull ($\alpha$ = shape, $\beta$ = scale)** | | | | | | | | | | | |
| $\alpha=0.5$, $\beta=1.0$ | 12.50 | 3.56 | 2.00 | 128.68 | 121.09 | 22.67 | 20.27 | 217.99 | 210.63 | 39.11 | 36.58 |
| $\alpha=1.5$, $\beta=0.5$ | 0.83 | 0.55 | 0.45 | 2.20 | 2.14 | 1.23 | 1.16 | 2.71 | 2.67 | 1.53 | 1.48 |
| $\alpha=1.5$, $\beta=1.5$ | 2.49 | 1.64 | 1.35 | 6.60 | 6.41 | 3.68 | 3.49 | 8.14 | 8.01 | 4.58 | 4.45 |
| **Inverse Gaussian ($\mu$ = mean, $\sigma$ = scale)** | | | | | | | | | | | |
| $\mu=0.5$, $\sigma=1.0$ | 1.05 | 0.62 | 0.50 | 3.99 | 3.83 | 1.63 | 1.52 | 5.42 | 5.30 | 2.21 | 2.12 |
| $\mu=2.5$, $\sigma=0.5$ | 15.09 | 4.42 | 2.50 | 124.01 | 117.48 | 28.00 | 25.06 | 187.97 | 182.76 | 47.37 | 44.41 |
| $\mu=2.0$, $\sigma=2.0$ | 5.49 | 2.69 | 2.00 | 26.80 | 25.61 | 9.13 | 8.42 | 37.81 | 36.92 | 13.26 | 12.64 |

**Table 2.** Log-likelihood values and parameter estimates for beta, Log-Lindley and generalized Log-Lindley models, with and without covariates

| | | $Y$ | | | $1-Y$ | | |
|---|---|---|---|---|---|---|---|
| **Models** | | **Beta($a$,$b$)** | **Log-Lindley,** $LL_1(\theta,\lambda)$ | **Generalized Log-Lindley,** $LL_p(\theta,\lambda)$ | **Beta($a$,$b$)** | **Log-Lindley,** $LL_1(\theta,\lambda)$ | **Generalized Log-Lindley,** $LL_p(\theta,\lambda)$ |
| (a) *Without covariates* | | | | | | | |
| **Log-likelihood** | | 76.1175 | 76.6042 | 83.2511 | 76.1175 | 69.0196 | 69.0195 |
| **Estimates** | | $a = 0.6125$<br>$b = 3.7979$ | $\lambda = 0.0343$<br>$\theta = 0.6907$ | $\lambda = 0.3824$<br>$\theta = 1.2694$<br>$p = 1.7819$ | $a = 3.7979$<br>$b = 0.6125$ | $\lambda = 4.16\times10^3$<br>$\theta = 5.9076$ | $\lambda = 2.66\times10^3$<br>$\theta = 5.9077$<br>$p = 1.0\times10^{-6}$ |
| **Models** | **Explanatory Variables** | **Beta** $(\mu,\phi)$ | **Log-Lindley,** $LL_1(\theta,\lambda)$ | **Generalized Log-Lindley,** $LL_p(\theta,\lambda)$ | **Beta** $(\mu,\phi)$ | **Log-Lindley,** $LL_1(\theta,\lambda)$ | **Generalized Log-Lindley,** $LL_p(\theta,\lambda)$ |
| (b) *With covariates and logit link for regression* | | | | | | | |
| **Log-likelihood** | | 87.7230 | 83.6526 | 98.2977 | 87.7230 | 96.7054 | 96.7624 |
| **Estimates** | (Intercept) | 1.8880 | 1.8422 | 2.7972* | -1.8880 | -2.7403 | -2.6807 |
| | ASSUME | -0.01214 | -0.005083 | $-5.03\times10^{-3}$ | 0.01214 | 0.03166 | 0.03191 |
| | CAP | 0.1780 | 0.05793 | 0.05734 | -0.1780 | -0.7586 | -0.7463 |
| | SIZELOG | -0.5115* | -0.2917* | -0.2894* | 0.5115* | 0.6962* | 0.6963 |
| | INDCOST | 1.2363* | 0.7122 | 0.6989* | -1.2363* | -3.6193* | -3.7226 |
| | CENTRAL | -0.01216 | -0.01971 | -0.01947 | 0.01216 | $5.53\times10^{-3}$ | $6.10\times10^{-3}$ |
| | SOPH | -0.003721 | $5.10\times10^{-4}$ | $1.875\times10^{-4}$ | 0.003721 | 0.03671 | 0.03558 |
| | | $\phi = 6.331$ | $\lambda = 0.0199$ | $\lambda = 0.3759$<br>$p = 3.3418$* | $\phi = 6.331$ | $\lambda = 108.58$* | $\lambda = 299.47$<br>$p = 1.0\times10^{-6}$ # |

* The estimated coefficient is significantly different from 0 with $p$-value < 0.02.

# Estimate approaches the boundary of parameter space and hence, inferences cannot be made for this set of coefficients.

**Table 3.** Log-likelihood values and parameter estimates for re-parameterized Log-Lindley and generalized Log-Lindley models with covariates

| Models | Explanatory Variables | $Y$ | | $1-Y$ | |
|---|---|---|---|---|---|
| | | **Log-Lindley,** $LL_1(\mu,\phi)$ | **Generalized Log-Lindley,** $LL_\gamma(\mu,\phi)$ | **Log-Lindley,** $LL_1(\mu,\phi)$ | **Generalized Log-Lindley,** $LL_\gamma(\mu,\phi)$ |
| **Log-Likelihood** | | 83.7575 | 91.9525 | 96.7689 | 96.8252 |
| **Estimates** | (Intercept) | 1.6767 | 0.2175 | -2.7200 | -2.7149 |
| | ASSUME | $-7.57\times10^{-3}$ | 0.004907 | 0.03208 | 0.03221 |
| | CAP | 0.08903 | -0.1335 | -0.7378 | -0.7508 |
| | SIZELOG | -0.4281* | -0.2644* | 0.7024* | 0.7049 |
| | INDCOST | 0.9687 | 0.4451 | -3.7854* | -3.7945 |
| | CENTRAL | -0.02318 | -0.07726 | $6.86\times10^{-3}$ | 0.002473 |
| | SOPH | $2.91\times10^{-4}$ | 0.005268 | 0.03593 | 0.03570 |
| | | $\phi = 0.03488$ | $\phi = 0.1755$ | $\phi = 67.212$* | $\phi = 65.0001$ |
| | | | $\gamma = 2.6735$* | | $\gamma = 1.0021^{\#}$ |

* The estimated coefficient is significantly different from 0 with $p$-value < 0.02.

[#] Estimate approaches the boundary of parameter space and hence, inferences cannot be made for this set of coefficients.

**Table 4.** Log-likelihood values and parameter estimates for beta, Log-Lindley and generalized Log-Lindley models, with and without covariates

| Models | Explanatory Variables | Beta($a$, $b$) | Log-Lindley, $LL_1(\theta,\lambda)$ | Generalised Log-Lindley, $LL_p(\theta,\lambda)$ |
|---|---|---|---|---|
| (a) *Without covariates* | | | | |
| **Log-Likelihood** | | 4197.526 | 3774.016 | 4333.53 |
| **Estimates** | | $a$ = 0.5462<br>$b$ = 27.965 | $\lambda$ = 0.0<br>$\theta$ = 0.3934 | $\lambda$ =1.8340<br>$\theta$ = 1.9236<br>$p$ = 8.0608 |
| (b) *With covariates and logit link for regression* | | | | |
| **Log-Likelihood** | | 4270.168 | 3794.483 | 4435.117 |
| **Estimates** | (Constant) | -5.0302* | -1.3478 | 0.3623* |
| | AGE | 0.01184* | 0.004766 | 0.004747* |
| | GENDER | 0.1170* | 0.05170 | 0.04789* |
| | insure | 0.3607* | 0.1509 | 0.1478* |
| | RACE(Asian) | -0.3493 | -0.1651 | -0.1984* |
| | RACE (Black) | -0.0130 | -0.02774 | -0.07190 |
| | RACE (Native) | 0.3063 | 0.1025 | 0.1535 |
| | RACE (White) | 0.02168 | -0.005097 | -0.04412 |
| | PHSTAT (Very Good) | 0.09941 | 0.04988 | 0.06882* |
| | PHSTAT (Good) | 0.1925* | 0.08882 | 0.1090* |
| | PHSTAT (Fair) | 0.4272* | 0.1877 | 0.2018* |
| | PHSTAT (Poor) | 0.7788* | 0.3069 | 0.3660* |
| | | $\phi$ = 32.453* | $\lambda = 0.0^{\#}$ | $\lambda = 1.97\times10^{-4}$<br>$p$ = 8.7734* |

* The estimated coefficient is significantly different from 0 with $p$-value < 0.02.

[#] Estimate approaches the boundary of parameter space and hence, inferences cannot be made for this set of coefficients.

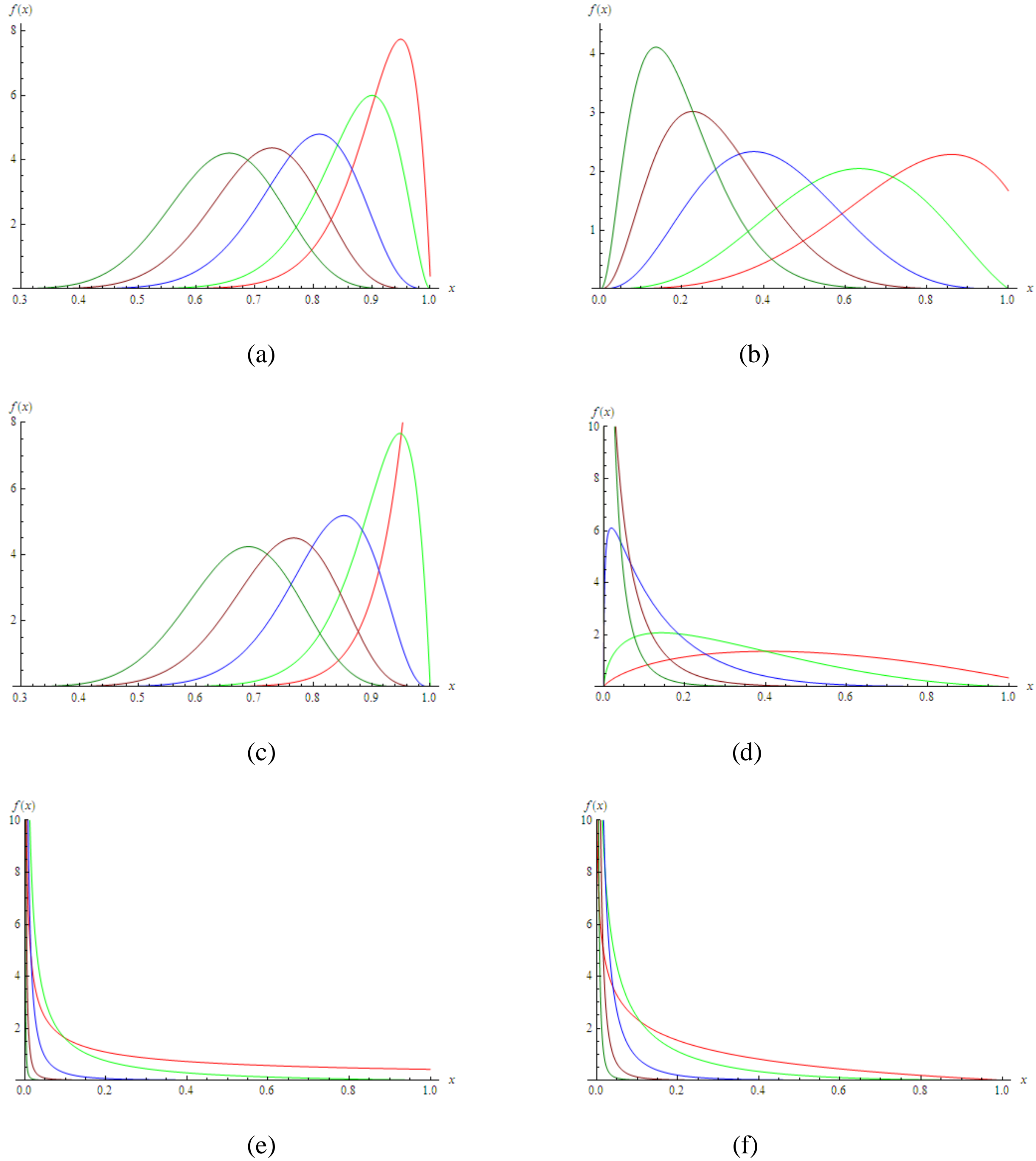

**Figure 1**. Pdf plots of $LL_p(\theta, \lambda)$ when $p = 0$ (red), 1 (green), 3(blue), 5 (brown), 7 (light green) for (a) $\theta = 20, \lambda = 0.001$ (b) $\theta = 5, \lambda = 0.1$ (c) $\theta = 20, \lambda = 5$ (d) $\theta = 2, \lambda = 0.1$ (e) $\theta = 0.5, \lambda = 10$ (f) $\theta = 0.9, \lambda = 0$

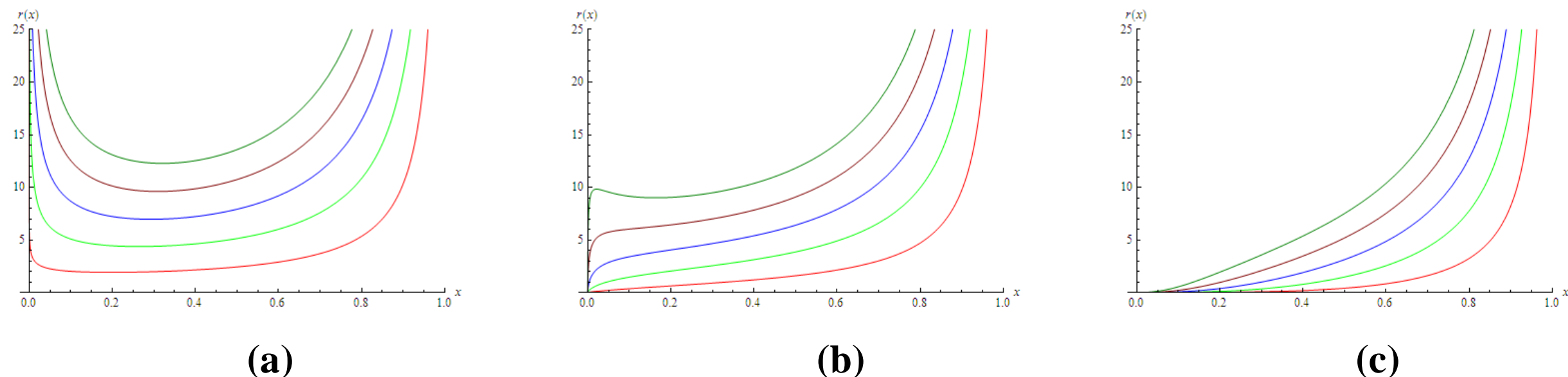

**Figure 2**. Hazard rate function plots of $LL_p(\theta,\lambda)$ for (a) $\theta=0.9,\lambda=2$ (b) $\theta=2,\lambda=2$ (c) $\theta=5,\lambda=2$ when $p=0$ (red), 1 (green), 3 (blue), 5 (brown), 7 (light green)